\documentclass[11pt]{amsart}

\usepackage{latexsym}
\usepackage{amssymb}
\usepackage{color}
\usepackage[dvips,frame,cmtip,arrow,matrix,line,graph]{xy}

\usepackage{amsmath}
\usepackage{amsthm}
\usepackage{amsfonts}
\usepackage{amscd}

\usepackage{epsfig}

\newtheorem{thm}{Theorem}[section]
\newtheorem{lem}[thm]{Lemma}
\newtheorem{prop}[thm]{Proposition}
\newtheorem{rem}[thm]{Remark}

\newtheorem{cor}[thm]{Corollary}

\newcommand{\vs}{\vspace{4mm}}

\newcommand{\CC}{\mathbb{C}}
\newcommand{\D}{\mathcal{D}}
\newcommand{\fD}{\mathfrak{D}}

\newcommand{\G}{\mathcal{G}}
\newcommand{\LL}{\mathbb{G}}
\newcommand{\M}{\mathcal{M}}
\newcommand{\N}{\mathbb{N}}

\newcommand{\ES}{\mathcal{S}}
\newcommand{\TT}{\mathfrak{T}}
\newcommand{\T}{\mathcal{T}}
\newcommand{\X}{\mathcal{X}}
\newcommand{\Z}{\mathbb{Z}}

\newcommand{\Simp}{\textsf{Simp}}

\newcommand{\Dif}{\textrm{Diff}^+}

\newcommand{\Ga}{\Gamma}

\newcommand{\Om}{\Omega}
\newcommand{\s}{\sigma}
\newcommand{\Si}{\Sigma}

\newcommand{\rar}{\longrightarrow}
\newcommand{\inc}{\hookrightarrow}
\newcommand{\sta}{\stackrel}

\newcommand{\x}{\times}

\newcommand{\sq}{\Box}
\newcommand{\cqfd}{\hfill $\square$}

\newcommand{\del}{\partial}

\newcommand{\Aut}{\operatorname{Aut}}

\newcommand{\B}{\operatorname{B}\!}

\newcommand{\Htpy}{\operatorname{Htpy}}
\newcommand{\Tel}{\operatorname{Tel}}
\newcommand{\colim}{\operatorname{colim}}

\newcounter{samcounter}

\begin{document}

\bibliographystyle{plain}

\title[From mapping class groups to automorphisms of free groups]{From mapping class groups to automorphism groups of free groups}

\author{Nathalie Wahl}
 \address{Department of Mathematics\\
               University of Aarhus\\
               Ny Munkegade 116\\
               8000 Aarhus\\
               DENMARK}
      \email{wahl@imf.au.dk}
      \thanks{Supported by a Marie Curie Fellowship of the European Community under contract number HPMF-CT-2002-01925}

\date{\today}

\begin{abstract}
We show that the natural map from the mapping class groups of surfaces to the automorphism groups of free groups,
induces an infinite loop map on the classifying spaces of the stable groups after plus construction.
The proof uses automorphisms of free groups with boundaries which play the role of mapping class groups of surfaces with several boundary components. 
\end{abstract}

\maketitle

\section{Introduction}

Both the stable mapping class group of surfaces $\Ga_\infty$ and the stable automorphism group of free groups $\Aut_\infty$ give rise to infinite loop spaces $\B\Ga_\infty^+$ and $\B\Aut_\infty^+$ when taking the plus-construction 
of their classifying spaces. By the work of Madsen and Weiss \cite{MadWei02}, $\B\Ga_\infty^+$
is now well understood, whereas $\B\Aut_\infty^+$ remains rather mysterious. 
In this paper, we relate these two spaces by showing that the natural map $\B\Ga_\infty^+\to \B\Aut_\infty^+$ is a map of infinite loop spaces. 
This means in particular that the map $H_*(\Ga_\infty)\to H_*(\Aut_\infty)$ respects the Dyer-Lashof algebra structure.
To prove this result, we introduce a new family of groups, the automorphism groups of free groups with boundary, which have the 
same stable homology as the automorphisms of free groups 
but enable us to define new operations.


Let $F_n$ be the free group on $n$ generators, and let $\Aut(F_n)$ be its automorphism group. Note that $\Aut(F_n) \cong \pi_0\Htpy_*(\vee_n S^1)$, 
the group of components of the pointed self homotopy equivalences of a wedge of $n$ circles. 
A result of Hatcher and Vogtmann \cite{HatVog-cerf} says that the inclusion map 
$\Aut(F_n)\to \Aut(F_{n+1})$ induces an isomorphism  $H_i(\Aut(F_n))\to H_i(\Aut(F_{n+1}))$ when  $i\le (n-3)/2$, and hence the homology of 
the stable automorphism group $\Aut_\infty:=\colim_{n\to\infty}\Aut(F_n)$ carries information about the homology of $\Aut(F_n)$ for $n$ large enough.

Let $S_{g,k}$ be a surface of genus $g$ with $k$ boundary components, and let $\Ga_{g,k}:=\pi_0\Dif(S_{g,k};\del)$ be its mapping class group, 
the group of components of the orientation preserving diffeomorphisms which fix the boundary pointwise. The maps 
 $\Ga_{g,k}\to \Ga_{g+1,k}$, induced by gluing a torus with two discs removed, and  
$\Ga_{g,k}\to \Ga_{g,k+1}$, induced by gluing a pair of pants, 
are homology isomorphisms in dimension $i\le (g-1)/2$ by Harer and Ivanov \cite{Har85,Iva90}. 
Let $\Ga_\infty:=\colim_{g\to\infty}\Ga_{g,1}$ denote the stable mapping class group.

There is a map  $f:\Ga_{g,1}\to \Aut(F_{2g})$ 
obtained by considering the action on the fundamental group of the surface $S_{g,1}$ since $\pi_1(S_{g,1})\cong F_{2g}$. 
The spaces $\coprod_{g\ge 0} \B\Ga_{g,1}$ and $\coprod_{n\ge 0} \B\Aut(F_n)$ have monoid structures induced by the pair of pants multiplication on surfaces
$S_{g,1}\x S_{h,1}\to S_{g+h,1}$ and by wedging circles $\vee_nS^1\x \vee_mS^1\to \vee_{n+m}S^1$ respectively. 
The map 
$\B f:\coprod_{g\ge 0} \B\Ga_{g,1}\to \coprod_{n\ge 0} \B\Aut(F_n)$ is a map of monoids, 
and hence we have a map of loop spaces $\Z\x \B\Ga_\infty^+ \to \Z\x \B\Aut_\infty^+$ on the group completion. 
The wedge product defines an infinite loop structure on $\Z \times \B\Aut_\infty^+$ as it defines a symmetric monoidal structure on $\coprod \Aut(F_n)$, 
thought of as a category.
The infinite loop structure on $\Z \times \B\Gamma_\infty^+$, discovered by Tillmann \cite{Til97}, is more complicated: Tillmann defines a 
cobordism 2-category $\ES$ which is symmetric monoidal under disjoint union and such that  $\Om\B\ES\simeq \Z\x\B\Ga_\infty^+$.
As $\B\ES$ is 
an infinite loop space, so is its loop space $\Z\x\B\Ga_\infty^+$. 
Our main result is the following.

\begin{thm}\label{main}
There is an infinite loop space structure on $\Z\x\B\Aut_\infty^+$ equivalent to the one induced by wedging circles and such that the 
 map $$\Z\x \B\Ga_\infty^+ \to \Z\x \B\Aut_\infty^+$$ induced by the action on the fundamental group is a map of infinite loop spaces. 
\end{thm}

To prove this theorem, we enlarge Tillmann's cobordism category $\ES$ by introducing an extra `graph-like' morphism.  
We obtain a new 2-category $\TT$ which contains $\ES$ as a subcategory and such that 
$\Om \B\TT\simeq \Z\x\B\Aut_\infty^+$. The theorem then follows from the fact that the inclusion $\ES\to \TT$ is a map of 
symmetric monoidal categories. 

The category $\TT$ is closely related to 
the {\em automorphism groups of free groups with boundary}, which we define now.

\begin{figure}[b]
\centering
\includegraphics[height=3cm]{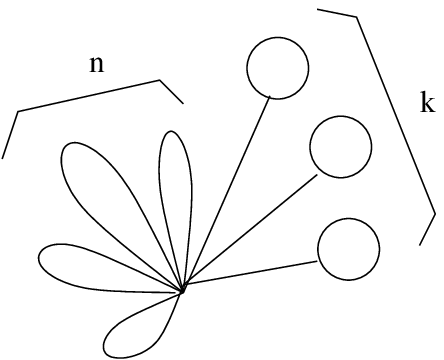} 
\caption{$G_{n,k}$}\label{Gnk}
\end{figure}

Let $G_{n,k}$ be the graph shown in Figure~\ref{Gnk} consisting of a wedge of $n$ circles together with $k$ extra circles joined by edges to the 
basepoint.
We call the $k$ circles disjoint from the basepoint {\em boundary circles}.
The automorphism group of free group with boundary $A_{n,k}$ is by definition $\pi_0\Htpy(G_{n,k};\del)$, 
the group of components of the space of homotopy 
equivalences of $G_{n,k}$ which fix the basepoint and the $k$ boundary circles pointwise. 
In particular, $A_{n,0}=\Aut(F_n)$. 
The group $A_{n,k}$ is an analogue of $\Ga_{g,k+1}$. In fact, 
$\Ga_{g,k+1}$ is a subgroup of $A_{2g,k}$. 
A description of $A_{n,k}$ as an extension of a subgroup of $\Aut(F_{n+k})$, as well as a presentation of the group, are given in joint work 
with Jensen \cite{JenWah04} 
and a description in terms of mapping class groups of certain 3-manifolds 
is given in joint work with Hatcher \cite{HatWah04}. 
This last description is used in \cite{HatWah04} to prove that the natural inclusions $A_{n,k}\to A_{n+1,k}$ and 
$A_{n,k}\to A_{n,k+1}$  
are homology isomorphisms in dimension $i\le (n-3)/3$.

The boundary circles of $G_{n,k}$ allow  to define 
new gluing operations between graphs and between graphs and surfaces. 
We use these operations to define the 2-category $\TT$, 
whose objects are the natural numbers, whose 1-morphisms are build out of graphs and surfaces by 
gluing and disjoint union, and whose 2-morphisms are homotopy equivalences fixing the boundary. 
The cobordism 2-category $\ES$ of \cite{Til97} is the subcategory of $\TT$ generated by the 1-morphisms build out of surfaces only. (See Section~\ref{category}.)

Theorem~\ref{LBT} says that $\Om \B\TT\simeq \Z\x \B\Aut_\infty^+$. 
The main ingredients of the proof 
are the homological stability of the automorphisms of free groups with boundary and 
a generalized group completion theorem. 
Our Theorem~\ref{right} then says that the infinite loop space structure on $\Z\x \B\Aut_\infty^+$  induced by the symmetric monoidal structure 
of $\TT$ is equivalent 
to the one previously known, induced by wedging circles.

Tillmann proved in \cite{Til99} that the map $\Z\x \B\Ga_\infty^+\to A(*)$ to Waldhausen's space $A(*)$ is a map of infinite loop spaces.
Our result says that 
this map factors as a map of infinite loop spaces: $\Z\x \B\Ga_\infty^+\to\Z\x \B\Aut_\infty^+\to A(*)$. The space $QS^0$ factors out of each of these 
three spaces (away from two in the case of the mapping class group) and the first map  on the 
$QS^0$ factor is multiplication by 2. 
The work of Dwyer-Weiss-Williams (\cite{DwyWeiWil}, see also \cite{Til99}) implies that the composite map $\Z\x\B\Ga_\infty^+ \to A(*)$ factors through 
$QS^0$. It is unknown whether $Z\x\B\Ga_\infty^+ \to \Z\x\B\Aut_\infty^+$ already factors through $QS^0$, or for that matter whether $\Z \times \B\Aut_\infty^+$
is $QS^0$. 
The rational homology of $\Aut_\infty$ is known to be trivial up to dimension 6 and is conjecturally trivial in all dimensions \cite{HatVog-rat}. 
On the other hand, the homology of the mapping class group is known with any field coefficients \cite{Gal04} and is rather rich.

\vs

The paper is organized as follows: 
We construct the 2-category $\TT$ in Section 2. In Section 3, we define a $\TT$-diagram   
which  is used in Section 4 
to show that $\Om\B\TT\simeq\Z\x\B\Aut_\infty^+$.  
In Section 5, we show the equivalence of the two infinite loop space structures on $\Z\x\B\Aut_\infty^+$. 
Section 6 gives an alternative definition of $\TT$ using punctured surfaces.

\subsection*{Acknowledgment}
I would like to thank the Aarhus group, and in particular Marcel B\"ocksted, S{\o}ren Galatius and Ib Madsen for
many helpful conversations. I would also like to thank  Allen Hatcher and Craig Jensen for their collaboration on 
the related papers \cite{HatWah04} and \cite{JenWah04}, 
and Benson Farb and Ulrike Tillmann for conversations at the early stages of this paper. 
Finally, I would like to thank the referee for suggesting a nice improvement of the paper.

\section{Cobordism category with graphs}\label{category}

The $(1+1)$-cobordism 2-category has objects 1-dimensional manifolds, 1-morphisms cobordisms between 
these manifolds and 2-morphisms diffeomorphisms of cobordisms restricting to the identity on the boundary.
We define here a 2-category $\TT$, modifying  Tillmann's model of the cobordism category in 
by adding an extra 1-morphism from the circle to itself, and replacing diffeomorphisms by homotopy 
equivalences.

\vs

\noindent{\bf Objects:}
The objects of $\TT$ are the natural numbers $n\in\N$, where $n$ can be thought of as a disjoint union of $n$ circles. 

\vs

\begin{figure}[ht]
\input{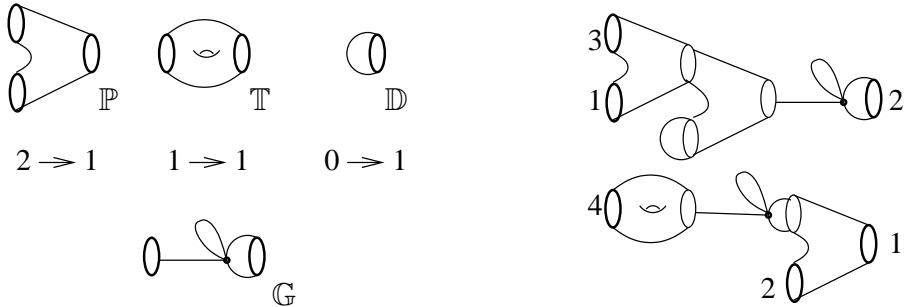} 
\caption{Pieces and example of a 1-morphism from 4 to 2}\label{pieces}
\end{figure}

\noindent{\bf 1-morphisms:}
To define the 1-morphisms of $\TT$, consider the following directed building blocks, called {\em pieces}:\\
- a pair of pants $\mathbb{P}$, a torus with two discs removed $\mathbb{T}$ and a disc $\mathbb{D}$ 
with respectively $2,1$ and 0 incoming boundary components and each with 1 outgoing boundary; \\
- a graph-like piece $\mathbb{G}\simeq  G_{1,1}\vee D^2$, with incoming boundary the boundary 
circle of $G_{1,1}$ and outgoing boundary  the boundary of the disc
(see Fig.~\ref{pieces}).

We fix  a parametrization $[0,2\pi[$ of the boundary circles, with
 0 at the end point of the attaching edge for the incoming boundary of $\LL$. 
We then allow the following gluing operation on the pieces: an incoming circle of one piece can be identified ---using the parametrization--- with the 
outgoing circle of another piece. The boundary of the glued object is defined to be the union of the boundary of the two pieces, minus the two identified
circles.

\vs

A 1-morphism in $\TT$ from $n$ to $m$ is  a couple $(T,\s)$, where $T$ is a 2-dimensional CW-complex obtained from the above 
pieces by gluing and disjoint union, with $n$ incoming and $m$ outgoing boundary components, and $\s$ 
is a labeling of these boundaries. (See Fig.~\ref{pieces} for an example.) 
One can moreover take disjoint union with copies of the circle, thought of as a morphism from 1 to 1. These circles are thus to be labeled on both sides. 
The $n!$  morphisms from $n$ to $n$ given by $n$ disjoint copies of the circle
correspond to the permutations of the labels, with the identity permutation representing the identity morphism on $n$. 

One should think of $(T,\s)$ as a combinatorial object, to which a topological space is associated.  

For simplicity, we will drop $\s$ from the notation in $(T,\s)$. 
Note that the building blocks of $\TT$ are defined in such way that a 1-morphism $T:n\to m$ has exactly $m$ connected components.

\vs

\noindent{\bf 2-morphisms:}
A 2-category is a category enriched over categories
and the 1-morphisms of $\TT$ are the objects of the  categories of morphisms $\TT(n,m)$. 
The 2-morphisms are the morphisms in these categories. Given two 1-morphisms $T$ and $T'$, objects of $\TT(n,m)$, the set of 
2-morphisms between $T$ and $T'$ is $\pi_0\Htpy(T,T';\del)$, the group of components of the 
space of homotopy equivalences from $T$  to $T'$ which fix the boundary, i.e. which map the boundary of $T$ to the boundary of $T'$ via the identity map 
according to the labels. 
By \cite[Prop. 0.19]{Hat02}, any homotopy equivalence $T\to T'$  which fixes the boundary is a homotopy equivalence relative to the boundary. 
(The pair $(T,\del T)$ is a CW-pair and hence satisfies the homotopy extension property, in which case the proposition applies.) 
In particular, every 2-morphism is invertible, but 
more importantly we will be able to glue homotopy equivalences along the boundaries.

\begin{rem}{\rm
We consider homotopy classes of homotopy equivalences relative to the boundary. By homotopy classes, we mean path components 
in $\Htpy(T,T';\del)$. Two elements $f$ and $g$ are in the same path component if and only if they are homotopy equivalent {\em relative to the boundary}. 
Note that here, it would not be equivalent to forget the ``relative to the boundary''!
}\end{rem}

Note that  if $T$ and $T'$ are permutations, the set of 2-morphisms is empty unless $T=T'$ in which case it is just the identity. 
If $T$ and $T'$ are surfaces, i.e. if they are build out of $\mathbb{P}$, $\mathbb{T}$ and $\mathbb{D}$ only, then 
$$\pi_0\Htpy(T,T';\del)\cong\pi_0\Dif(T,T';\del)$$ 
(see \cite[Sect. 2]{JenWah04}). We describe in Lemma~\ref{new} below the 1-morphisms with automorphism group given by automorphisms of free groups
with boundaries.

\vs

\noindent
{\bf Composition:}
The composition of 1-morphisms is defined by gluing according to the labels. 
We will denote this composition by $T_1\sq T_2:n\to p$ for $T_1:n\to m$ and $T_2:m\to p$.

For $f_1:T_1\to T_1'$ and $f_2: T_2\to T_2'$, define the horizontal composition $f_1\sq f_2: T_1\sq T_2\to T_1'\sq T_2'$ 
by applying $f_1$ to $T_1$ and $f_2$ to $T_2$. This defines a homotopy equivalence as $f_1$ and $f_2$ are homotopy equivalences relative to the 
boundary.

Finally, vertical composition of 2-morphisms is given by the composition of homotopy equivalences.

All compositions are associative and this defines a 2-category. 
We write all compositions in $\TT$ in the left-to-right order.

\vs

We have seen that surface 1-morphisms in $\TT$ have automorphism group given by mapping class groups. 
We now describe the 1-morphisms with automorphism group given by automorphisms of free groups with boundary. 
They play an important role in determining the homotopy type of $\TT$. 

Let $G_{n,k}$ be the graph given in Fig.~\ref{Gnk} and consider $G_{n,k}\vee D^2$, where  the basepoint of the disc is its center. 
The boundary of $G_{n,k}\vee D^2$ consists of $k+1$ circles: the $k$ boundary circles of $G_{n,k}$ together with the boundary circle of the disc $D^2$. 
The graph $G_{0,0}$ is just a point so $G_{0,0}\vee D^2\cong \mathbb{D}$. Note also that $G_{1,1}\vee D^2=\mathbb{G}$.

Recall that $A_{n,k}=\pi_0\Htpy(G_{n,k};\del)$, where $\del G_{n,k}$ is the basepoint of $G_{n,k}$ union its $k$ boundary circles. 

\begin{figure}[ht]
\input{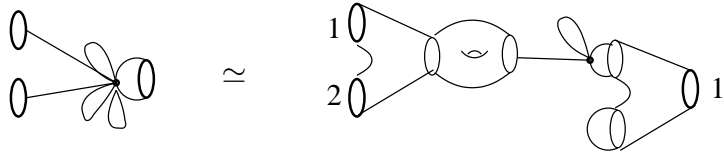} 
\caption{$G_{3,2}\vee D^2$ homotopic to an object of $\TT_D(2,1)$}\label{discs}
\end{figure}

We call {\em surface component} of an object $T$ of $\TT(k,1)$ any connected component of surface pieces $\mathbb{P,T,D}$ and the disc 
of $\LL$
glued together in $T$. For example, the object of $\TT(2,1)$ in Fig.~\ref{discs} has two surface components: a torus with three holes and a disc.
The {\em outgoing surface component} of $T$ is the surface component whose outgoing boundary is the outgoing boundary of $T$. 
In the example, it is a disc.

\begin{lem}\label{new}
An object $T$ of $\TT(k,1)$ is homotopic relative to the boundary to $G_{n,k}\vee D^2$ for some $n$,  if and only if
the outgoing surface component of $T$ is a disc.  
Moreover, if this is the case, we have $\ \pi_0\Htpy(T;\del)\cong A_{n,k}$.
\end{lem}

Let $\TT_D(k,1)$ be the full subcategory of $\TT(k,1)$ generated by the objects with a disc as outgoing surface component, 
and let $\TT_S(k,1)$ be the full subcategory of $\TT(k,1)$ generated 
by the other objects, i.e. the objects  with outgoing surface component of higher genus or with more boundary components, or the 
object represented by the circle in $\TT(1,1)$. 
By the lemma, there can be no morphisms in $\TT(k,1)$ between the objects of $\TT_D(k,1)$ and the objects of $\TT_S(k,1)$. Hence we have

\begin{cor}\label{split}
$$\TT(k,1)=\TT_D(k,1)\coprod \TT_S(k,1)\ \ \  {\rm and}\ \ \ \B\TT_D(k,1)\simeq \coprod_{n\ge \varepsilon} \B A_{n,k},$$
where $\varepsilon=0$ when $k=0$ and $\varepsilon=1$ otherwise.
\end{cor}

\noindent
{\em Proof of Lemma \ref{new}.}
Every surface component of an object $T$ has exactly one outgoing boundary by construction.  
Suppose that $T$ is an object of $\TT(k,1)$ with outgoing surface component a disc.  
Then the outgoing boundary of any other surface component in $T$ will be free, i.e.~not fixed by the homotopy equivalences 
---as it cannot be the outgoing boundary of $T$.
One can thus homotope each of these surface components $S$ to a graph $G_{2m,l}$ lying in $S$, with basepoint $0$ of the outgoing boundary of $S$, its $l$ 
boundary circles on the $l$ 
incoming boundaries of $S$ and where $m$ is the  genus of $S$, and this can be done while keeping the boundary of $T$ fixed. Such a 
homotopy can easily be constructed by considering the surface component as a polygon with appropriate identifications of the edges and with 
$l+1$ discs removed. 


Suppose on the other hand that $T$ does not end with a disc component and let $S$ be the outgoing surface component 
of $T$. 
If $S$ has incoming boundary components, then they have to be incoming boundaries of $T$ because 
the only non-surface piece one can glue to $S$ is $\LL$ and this closes the incoming boundary it is glued to. Hence  
all the boundary components of $S$ are boundary components of $T$ and  
 $S$ cannot be homotoped to 
a graph while fixing the boundary of $T$. (Recall that $\pi_0\Htpy(S;\del)\cong\pi_0\Dif(S;\del)$ for a surface $S$.)


Finally, we want to see that if $T$ is homotopic to $G_{n,k}\vee D^2$ relative to the boundary, then $\pi_0\Htpy(T;\del)\cong A_{n,k}$. 
It is equivalent to show that $\pi_0\Htpy(G_{n,k}\vee D^2;\del)\cong A_{n,k}$. Let $X:=G_{n,k}\vee D^2$. 

Note first that $A_{n,k}\cong \pi_0\Htpy(X;\tilde\del)$, where $\tilde\del X$ consists of the $k$ boundary circles 
of $G_{n,k}$ together with the basepoint $*=0\in \del D^2$.  There is an inclusion 
$$\Htpy(X;\del)\inc \Htpy(X;\tilde\del)$$
and we want to show that it induces an isomorphism on $\pi_0$. 

\noindent
{\em Surjectivity:} Let $f\in \Htpy(X;\tilde\del)$. Then $f$ maps $\del D^2$ to a trivial loop in $X$. 
It is then easy to homotope $f$ relative to $\tilde\del X$ to a map fixing $\del D^2$. (One can also use the homotopy extension property of $(X,\del X)$.)

\noindent
{\em Injectivity:} Suppose $f,g\in \Htpy(X;\del)$ are homotopic relative to $\tilde \del$ by a homotopy $H:X\x I\to X$. Consider the restriction of $H$
$$H_{\del D^2}: \del D^2\x I\to X.$$ 
Let  $\hat H_{\del D^2}$ be a lift of $H_{\del D^2}$ to the universal cover of $X$, which is contractible. As $H_{\del D^2}(0\x I)=*$ and 
$H_{\del D^2}(\del D^2\x 0)=H_{\del D^2}(\del D^2\x 1)=\del D^2$, we have 
$\hat H_{\del D^2}(\del D^2\x 0)=\hat H_{\del D^2}(\del D^2\x 1)$ is a fixed lift $\hat{\del D^2}$ of $\del D^2$. 
Hence $\hat H_{\del D^2}$ is homotopic to the constant map 
$i_{\del D^2}: \del D^2\x I\to \hat{\del D^2}$  relative to $\del D^2\x \{0,1\}$. 
Let $h:(\del D^2\x I)\x I\to X$ be the projection of this last homotopy and extend it to  
$h_\del:((\del X\x I)\cup (X\x\{0,1\})) \x I\to X$ using the constant map on $\del X\backslash \del D^2$ and 
$f$ and $g$ on $(X\x 0)\x I$ and $(X\x 1)\x I$. 
Hence $h_\del$ is a homotopy from the restriction $H_{\del}$ of $H$ to $(\del X\x I)\cup (X\x\{0,1\})$ 
to the map $i_{\del X}\cup f\cup g:(\del X\x I)\cup (X\x\{0,1\}) \to X$.  
Now $(X\x I ,(\del X\x I)\cup (X\x\{0,1\}))$ has the homotopy extension property as it is a CW-pair, so we get a homotopy 
${\bf h}:(X\x I)\x I \to X$ from $H$ to a homotopy $H': X\x I\to X$ with $H'(X\x 0)=f(X)$, $H'(X\x 1)=g(X)$ and 
$H'(\del X\x I)$ is the constant map on $\del X$. Hence $f$ and $g$ are homotopic relative to $\del X$. 
\cqfd

\vs

By a $\Delta$-category, we mean a category enriched over simplicial sets.

Let $\T$ denote the $\Delta$-category obtained from $\TT$ by taking the nerve of the categories of morphisms, that is $\T$ has the same objects 
as $\TT$, and $\T(n,m):=N_\bullet\TT(n,m)$.

\begin{prop}
The $\Delta$-category $\T$ is symmetric monoidal under disjoint union and its classifying space $\B\T$ is an infinite loop space.
\end{prop}

\begin{proof}
Taking disjoint union $T_1\sqcup T_2$ and shifting the labels of $T_2$ induces a monoidal structure on $\T$. 
The symmetries are then given by the block permutations $n+m\to m+n$ 
in $\T(n+m,n+m)$. 

Now $\B\T$ is connected as there is a morphism from any object to 1. It follows that $\B\T$ is an infinite loop space \cite{May74,Seg74}. 
\end{proof}

The cobordism $\Delta$-category $\ES$ defined in \cite{Til99} is the subcategory of $\T$ generated by the surface pieces $\mathbb{P,T}$ and 
$\mathbb{D}$, i.e. it has the same objects 
as $\T$ and the morphism space $\ES(n,m)$ is the nerve of the full subcategory of $\TT(n,m)$ generated by the surface objects.
(To be precise, the above model for $\ES$ is a mix of \cite{Til97} and \cite{Til99}: we use the combinatorial description of \cite{Til99} but 
work with mapping class groups as in \cite{Til97}, instead of spaces of diffeomorphisms. This is equivalent because  
the components of the diffeomorphism group of a surface are contractible 
when the genus is large enough \cite{EarEel,EarSch}.)

The $\Delta$-category $\ES$ is symmetric monoidal under disjoint union and 
$$\Om \B\ES\simeq \Z\x\B\Ga_\infty^+$$
\cite[Thm. 3.1]{Til97}.

As the inclusion respects the symmetric monoidal structure, we have the following proposition:

\begin{prop}\label{infinit1}
The inclusion of categories $\ES\to \T$ induces a map of infinite loop spaces $\B\ES\to \B\T$.
\end{prop}

Theorem~\ref{LBT} says that $\Om\B\T\simeq\Z\x\B\Aut_\infty^+$. 
The map of infinite loop spaces $\Om\B\ES\to \Om\B\T$ is then the map announced in Theorem~\ref{main}.

\section{A $\T$-diagram}\label{T-diagram}

In this section, we work with simplicial sets and bisimplicial sets---the latter coming from homotopy colimits of functors to the category of simplicial sets. 
We use the following fact about bisimplicial sets:
to any bisimplicial set $X_{\bullet,\bullet}$, one can associate the {\em diagonal simplicial set} $dX_\bullet$ by taking $dX_p=X_{p,p}$. 
Let $f:X_{\bullet,\bullet}\to Y_{\bullet,\bullet}$ be a bisimplicial map. If   
 $f_{p,\bullet}$ 
is a homotopy equivalence on the vertical simplicial sets 
$X_{p,\bullet}$ and $Y_{p,\bullet}$ for each $p\ge 0$, then the induced map 
$f:dX_\bullet\to dY_\bullet$ is also a homotopy equivalence \cite[IV Lem.~2.6]{GoeJar99}. 

When we consider a bisimplicial set as simplicial set, we mean its associated diagonal simplicial set. 

\vs

Given a small $\Delta$-category $\M$, an {\em $\M$-diagram} is a functor $\X:\M^{op}\to \Simp$, from the opposite category of $\M$ to 
the category of simplicial sets. It is thus 
 a simplicial set $\mathcal{X}=\coprod_i\mathcal{X}(i)$ where $i$ 
runs over the objects of $\M$, with a simplicial action
$$\M(i,j)\times\mathcal{X}(j) \rar \mathcal{X}(i)$$
for all objects $i,j$ of $\M$. 
An example of an $\M$-diagram is given by taking $\X(i)=\M(i,i_0)$ for a fixed object $i_0$ of $\M$.

We denote by $\T_1$ the $\T$-diagram with $\T_1(k):=\T(k,1)$ and  
by $\T_\infty$ the $\T$-diagram with 
$$\T_\infty(k)=\Tel(\T(k,1)\stackrel{\sq\LL}{\rar}\T(k,1)\stackrel{\sq\LL}{\rar}\dots),$$ 
with the maps in the telescope defined 
by composition with $\LL\in\T(1,1)$, and the $\T$-diagram structure, as for $\T_1$, obtained by pre-composition in $\T$.

Let $A_{\infty,k}:= \colim(A_{n,k}\to A_{n+1,k}\to\dots)$, where the map $A_{n,k}\to A_{n+1,k}$ is given by wedging with an $S^1$ and extending 
the homotopies via the identity.  

\begin{thm}\label{type}
$\T_\infty(k) \simeq \Z\x \B A_{\infty,k}$
\end{thm}

\begin{proof}
Recall from Corollary~\ref{split} that $\TT(k,1)= \TT_D(k,1) \sqcup \TT_S(k,1)$, where  
 $\TT_D(k,1)$ is the full subcategory of $\TT(k,1)$ generated by objects with a disc as outgoing surface component. 
Let $\T_D(k,1):=N_\bullet \TT_D(k,1)$ and $\T_S(k,1):=N_\bullet \TT_S(k,1)$ denote their nerves.

Let $\T_{D,\infty}(k)\subset \T_\infty(k)$ be the restriction 
of the telescope $\T_\infty(k)$ to the space $\T_D(k,1)$. 
As $\T_D(k,1)\simeq \coprod_{n\ge \varepsilon}\B A_{n,k}$ (Corollary~\ref{split}), we have $$\T_{D,\infty}(k)\simeq \Z\x \B A_{\infty,k}.$$ 
The theorem then follows from the fact that there is a retraction $$r:\T_\infty(k)\stackrel{\simeq}{\rar} \T_{D,\infty}(k)$$  
as the maps defining the telescope take $\T_S(k,1)$ to $\T_D(k,1)$: 
$$\T_\infty(k)=\Tel\Big(\ \ \xymatrix{\T_S(k,1)\ar@{.}[d]|-{\coprod}\ar[dr]^{\sq\LL} & \T_S(k,1)\ar@{.}[d]|-{\coprod}\ar[dr]^{\sq\LL} 
       & \T_S(k,1)\ar@{.}[d]|-{\coprod} \ar[dr] & \dots \ar@{.}[d]\\
\T_D(k,1)\ar[r]^{\sq\LL} & \T_D(k,1)\ar[r]^{\sq\LL} & \T_D(k,1) \ar[r] & \dots 
} \ \ \Big)$$
Explicitly, $r$ is a bisimplicial map defined vertically by 
$r_p:(\T_\infty(k))_p=\coprod_{n_0\le\dots\le n_p\in\N} (\T_S(k,1)\sqcup \T_D(k,1)) \to (\T_{D,\infty}(k))_p=\coprod_{n_0\le\dots\le n_p\in\N}\T_D(k,1)$
which maps $\T_S(k,1)\x (n_0,\dots,n_0,n_1,\dots,n_{p-i})$ with $n_0<n_1$ to $T_D(k,1)\x (n_0+1,\dots,n_0+1,n_1,\dots,n_{p-i})$ using the map 
$\sq\LL$ 
and maps $\T_D(k,1)\x (n_0,\dots,n_{p})$ to itself via the identity.
\end{proof}

For an $\M$-diagram $\X$, we can consider the Borel construction $E_\M\X$ of the action, 
that is the homotopy colimit of the functor $\X$. 
So $E_\M\X$ is the bisimplicial set whose simplicial set of 
$p$-simplices is given by 
$$(E_\M\X)_p:=\coprod_{i_0,\dots,i_p\in Ob(\M)}\M(i_0,i_1)\x\dots\x\M(i_{p-1},i_p)\x\X(i_p)$$
with boundary maps $d_0$ by dropping, $d_1,\dots,d_{p-1}$ by composition in $\M$ and $d_p$ using the $\M$-diagram structure.

\begin{lem}\label{lemcontract}
$E_\T\T_\infty$ is contractible.
\end{lem}

\begin{proof}
By \cite[Lem.~3.3]{Til97}, $E_\T\T_1$ is contractible, where $\T_1$ is the $\T$-diagram with $\T_1(n)=\T(n,1)$. Then 
$E_\T\T_\infty\simeq \Tel(E_\T\T_1\to E_\T\T_1\to\dots)$ is a telescope of contractible spaces and thus is itself contractible. 
\end{proof}

\section{Homotopy type of $\TT$}\label{htpytype}

We use the $\T$-diagram $\T_\infty$ and 
the homological stability of the automorphisms of free groups with boundary $A_{n,k}$ \cite{HatWah04} to prove the following theorem:

\begin{thm}\label{LBT}
$\Om \B\T \simeq \Z\x \B\Aut_\infty^+$
\end{thm}

\begin{lem}\label{homo}
The vertices of the simplicial sets of morphisms in $\T$ act by homology isomorphisms on $\T_\infty$. 
\end{lem}

\begin{proof}
The vertices of $\T(n,m)$ are the objects of the category $\TT(n,m)$, so they are 2-dimensional CW-complexes built out of the 
pieces $\mathbb{P},\mathbb{T},\mathbb{D}$ and $\LL$ by gluing and disjoint union. It is enough to show that each of these building blocks 
acts by homology isomorphisms. An element of $\T_\infty(k)\simeq\Z\x\B A_{\infty,k}$ is a colimit of morphisms from $k$ to $1$ and $\T$ acts 
by precomposing, that is by gluing on the $k$ incoming boundary circles. 
Acting one piece at a time means precomposing with that piece disjoint union with $k-1$ circles.
We show that the corresponding operations on the groups $A_{n,k}$ induce homology isomorphisms stably.

By \cite{HatWah04}, we know that the maps 
$$A_{n,k}\stackrel{G_{1,0}}{\rar} A_{n+1,k}$$ 
$$A_{n,k}\stackrel{G_{0,1}}{\rar} A_{n,k+1}$$
obtained by wedging a circle (that is the graph $G_{1,0}$)
and a boundary circle (that is the graph $G_{0,1}$), induce homology isomorphisms on the stable groups:
$(G_{1,0})_*:H_*(A_{\infty,k})\stackrel{\cong}{\rar} H_*(A_{\infty,k})$ and 
$(G_{0,1})_*:H_*(A_{\infty,k})\stackrel{\cong}{\rar} H_*(A_{\infty,k+1}).$

This implies that gluing a disc $\mathbb{D}$ also induces 
a homology isomorphism stably since the 
composition 
$$A_{n,k-1}\stackrel{G_{0,1}}{\rar} A_{n,k} \stackrel{\mathbb{D}}{\rar} A_{n,k-1},$$
gluing $\mathbb{D}$ on the boundary circle of the added $G_{0,1}$, is the identity. 

Thus gluing a pair of pants $\mathbb{P}$ induces a homology isomorphisms stably as 
the composition with gluing discs on the incoming boundaries of $\mathbb{P}$ 
is the same as gluing a disc $\mathbb{D}$ where $\mathbb{P}$ was glued.

Gluing a torus $\mathbb{T}$  also induces a homology isomorphism stably as 
the composition 
$$A_{n,k-1}\stackrel{G_{0,1}}{\rar} A_{n,k} \stackrel{\mathbb{T}}{\rar} A_{n+2,k}$$
gluing $\mathbb{T}$ on the boundary circle of $G_{0,1}$, 
is the same as wedging the graph $G_{2,1}$. 
Indeed, the two  maps 
add two  extra circles and a boundary circle and prolong the homotopy equivalences via the identity on the added pieces.
Similarly, gluing $\LL$ induces a homology isomorphism stably as the composition with wedging $G_{0,1}$ is the same as 
wedging $G_{1,1}$. 
\end{proof}

Note that 
the three maps $A_{n,k}\to A_{n+2,k}$ induced by gluing $\mathbb{T}$, gluing $\LL\sq\LL$ or wedging $G_{2,0}$, are all different. 
However, by freeing a boundary circle (which induces a homology isomorphism stably)
one can show that the compositions $F\circ \mathbb{T}$ and $F\circ G_{0,2}:A_{n,k}\to A_{n+2,k}\to A_{n+3,k-1}$, 
 gluing a torus or wedging two circles and then freeing the appropriate boundary circle, are conjugate, 
and hence induce the same map on homology. 
This is another way to see that gluing $\mathbb{T}$ induces a homology isomorphism stably.

\vs

\noindent {\em Proof of Theorem~\ref{LBT}.}
Consider the square 
$$\xymatrix{\T_\infty(k)\ar[d]\ar[r] & E_\T\T_\infty \ar[d]^p \\
             k \ar[r] & \B\T.}$$
By the generalized group completion theorem \cite[Thm.~3.2]{Til97} (see also \cite{Moe89,PitSch04}) and using Lemma~\ref{homo}, this square is homology 
cartesian, which means that the fiber of $p$ 
at any vertex is homology equivalent to the homotopy fiber. 
As $E_\T\T_\infty$ is contractible (Lemma~\ref{lemcontract}), the 
homotopy fiber is $\Om \B\T$. Considering the fiber of $p$ at 0, we thus have a homology equivalence $\T_\infty(0)\to \Om \B\T$, and hence again a homology 
equivalence after plus-construction with respect to any perfect subgroup of $\pi_1\T_\infty(0)$. 
Now $\T_\infty(0)\simeq\Z\x \B\Aut_\infty$ by Thm.~\ref{type}
and 
the commutator subgroup $[\pi_1\T_\infty(0),\pi_1\T_\infty(0)]$ is perfect \cite[Rem.~2]{McDSeg}.
So we can plus-construct with respect to this subgroup and 
we get a homotopy equivalence by Whitehead's theorem for simple spaces \cite[Ex.~4.2]{Dro71}.
\cqfd

\section{Two infinite loop space structures on $\Z\x\B\Aut_\infty^+$}\label{two}

Consider the category $\TT_D(0,1)$. Its classifying space is homotopy equivalent to $\coprod_{n\ge 0}\B\Aut(F_n)$ by Corollary~\ref{split}. 
We claim that the pair of pants defines a symmetric monoidal structure (up to homotopy) on $\TT_D(0,1)$ equivalent to the 
symmetric monoidal structure induced by wedging circles on $\coprod_{n\ge 0}\Aut(F_n)$, thought of as a category with objects 
$\{\vee_nS^1|n\in\N\}$ and morphism sets $\pi_0\Htpy_*(\vee_n S^1,\vee_mS^1)$. 
To be able to remove the `up to homotopy' in the above claim, but also for the  comparison to 
the infinite loop space structure of $\Om\B\T$, we need to work with a quotient of the 2-category $\TT$. 

To make the pair of pants multiplication unital, we need to collapse the 1-morphisms $\mathbb{D}\sq_i\mathbb{P}$ to the circle, where $i=1,2$ 
represents the two possible gluings of $\mathbb{D}$ and $\mathbb{P}$. 
To make it associative, we need to identify $\mathbb{P}\sq_1\mathbb{P}$ and $\mathbb{P}\sq_2\mathbb{P}$. 
We define a 2-category $\TT^r$ whose objects are the natural numbers as in $\TT$, and whose 1-morphisms are 
the 1-morphisms of $\TT$
with no occurrence of the sequences $\mathbb{D}\sq_i\mathbb{P}$ for $i=1,2$, or $\mathbb{P}\sq_2\mathbb{P}$. The 2-morphisms of $\TT^r$ are 
as in $\TT$ except that the circle is thought of as a small cylinder, having thus automorphism group $\Z$. 
As shown in  \cite[Sec.~3.1.1]{Wah04}, for each object $T$ of $\TT(n,m)$ there exists a canonical quotient object $T^r$ of
$\TT^r(n,m)$ and a homotopy equivalence $T\to T^r$ which can be used to define composition in   
 $\TT^r$. (The argument in \cite{Wah04} extends immediately from $\ES$ to $\T$ as the quotient construction only affects the pieces $\mathbb{P,T}$ and
$\mathbb{D}$.)
We have 
$\Om\B\T^r\simeq\Z\x\B\Aut_\infty^+$ by the same argument as for $\T$. We chose to work with $\TT$ rather than $\TT^r$ in the first part of the 
paper to stay as simple as possible.

The pair of pants induces a monoidal structure on $\TT^r_D(0,1)$, where $\TT^r_D(0,1)$ is defined in an analogous way to $\TT_D(0,1)$. 
(Note that the quotient construction only affects the disc $\mathbb{D}$ and not the disc component in $\LL$.)  
The twist on the pair of pants defines a 
symmetry as it squares to a Dehn twist along the outgoing boundary of the pair of pants, which 
is homotopically trivial once the incoming boundaries of the 
pair of pants are closed by discs.

\begin{thm}\label{right}
The equivalence $\Om \B\T\simeq \Z\x \B\Aut_\infty^+$ is an equivalence of infinite loop spaces, where the left infinite loop structure is 
induced by disjoint union in $\T$ and the right one by wedging circles in $\coprod_{n\in\N}\Aut(F_n)$.
\end{thm}

\begin{proof}
The proof of this theorem is totally analogous to the proof of the main result in \cite{Wah04}, which says that the infinite loop space 
structure on $\Z\x \B\Ga_\infty^+$ coming from $\Om \B\ES$ is equivalent to the one coming from the $\M$-algebra structure of 
$\ES(0,1)\simeq\coprod_{g\ge 0}\B\Ga_{g,1}$, where $\M$ is the 
mapping class group operad defined in \cite{Til00}. 
We will only sketch the main steps for the case of interest here. 

We have an equivalence of infinite loop spaces $\Om\B\T\sta{\simeq}{\rar}\Om\B\T^r$, where $\T^r$ denotes the $\Delta$-category associated to $\TT^r$, 
and we have an equivalence of symmetric monoidal categories 
$\TT^r_\D(0,1)\sta{\simeq}{\rar} \coprod_{n\in\N}\Aut(F_n)$. 
We give here an equivalence between the spectra of deloops of $\Om\B\,(\T^r_D(0,1))$ and of $\Om\B\T^r$.

Let $\Si$ be the standard $E_\infty$ operad with $\Si(k)=E\Si_k$ 
for $\Si_k$ the symmetric group. (This operad is denoted $\Ga^+$ in \cite{BarEcc1,BarEcc2} and $\Ga$ in \cite{Wah04}.)
As $\TT^r_D(0,1)$ and $\T^r$ are symmetric monoidal categories,  the operad $\Si$ acts on their classifying spaces $\T^r_D(0,1)$ and $\B\T^r$. 

We want to 
relate $\Om \B\T^r$ and $\T^r_D(0,1)$, which is a subspace of $\T^r(0,1)$. 
To a morphism from 0 to 1 in $\T^r$ corresponds a 1-simplex with boundary points 0 and 
1 in $\B\T^r$, by definition of the nerve.
There is thus a  natural map 
 $$\phi:\T^r_D(0,1)\to \Om \B\T^r$$ taking the path from 0 to 1 in $\B\T^r$ defined by the element of $\T^r_D(0,1)$ 
and going back to 0 along the path defined by the disc, 
also thought of as a morphism from 0 to 1. The crucial observation is that this map respects the multiplication up to homotopy, 
as show in Figure \ref{htpy}.  
\begin{figure}[ht]
\includegraphics[width=9cm]{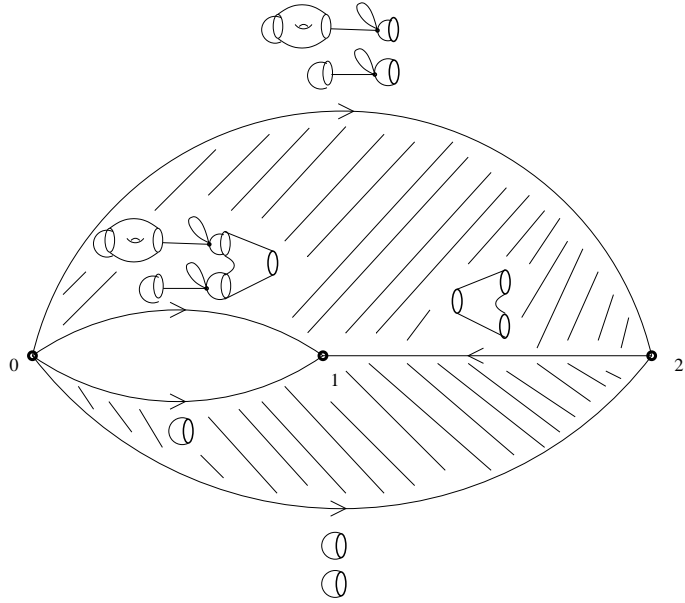}
\caption{Homotopy}\label{htpy}
\end{figure}
The figure shows the loop obtained by multiplying  $\mathbb{D}\sq\mathbb{T}\sq\LL$ and $\mathbb{D}\sq\LL$ in $\T^r_D(0,1)$ and then map to $\Om \B\T^r$, 
producing 
a loop following $(\mathbb{D}\sq\mathbb{T}\sq\LL\sqcup \mathbb{D}\sq\LL)\sq \mathbb{P}$ from 0 to 1, and then following the disc-morphism back to 0. 
The other loop, going from 0 to 2 and back to 0, 
is the one obtained by mapping first $\mathbb{D}\sq\mathbb{T}\sq\LL$ and $\mathbb{D}\sq\LL$ to $\Om \B\T^r$ and then multiplying, 
that is taking loop on disjoint union. 
The two loops are homotopic in $\B\T^r$ because the operad $\Si$ acts on $\T^r_D(0,1)$ by taking disjoint union and composing with $\mathbb{P}$ 
as a morphism in $\T^r$. 
This means that the top triangle in the  figure commutes in $\T^r$, and hence defines a 2-simplex in $\B\T^r$. On the other hand, we have defined 
composition in $\T^r$ so that the bottom triangle also commutes. (This is a place where we need the disc to be a strict unit.)

This argument extends to give all necessary higher homotopies, parametri\-zed by the elements of $\Si$. (Considering higher multiplications 
will give higher simplices in $\B\T^r$). 

More concretely, we consider Barratt and Eccles' model for the spectra associated to the action of $\Si$ on $\T^r_D(0,1)$ and on $\B\T^r$. 
These are realization of simplicial spaces, where the $i$th space of the spectrum has $q$-simplices given by 
$\G \Si (S^i\wedge \Si^q(\T^r_D(0,1)))$ and $\G\Si (S^i\wedge \Si^q(\B\T^r)))$ respectively, 
where $\G$ is the group completion functor on free simplicial monoids, 
$\Si$ denotes now the monad associated to the operad $\Si$, and $\Si^q$ means iterating the monad $q$ times. The simplicial structure is given by the 
monad multiplication for $d_0,\dots,d_{q-1}$ (using the assembly map $S^1\wedge \Si(X)\to \Si(S^1\wedge X)$ for $d_0$)
and by the $\Si$-algebra structure of $\T^r_D(0,1)$ and $\B\T^r$ for the last boundary map. 
The spectrum structure comes from the equivalence $\G\Si(X)\stackrel{\simeq}{\rar}\Om\G\Si(S^1\wedge X)$ \cite{BarEcc1,BarEcc2}.

The map $\phi$ defined above gives a map on the level of $q$-simplices from the $i$th deloop to the $(i-1)$st deloop:
$$\xymatrix{\G\Si (S^i\wedge \Si^q(\T^r_D(0,1)))\ar[r]^{f^i_q} \ar[d] &  \G\Si (S^{i-1}\wedge \Si^q(\B\T^r))\\
\G\Si (S^{i-1}\wedge \Si^q(S^1\wedge\T^r_D(0,1))). \ar[ur]^\phi
}$$ 
The maps $\{f^i_q\}_{q\ge 0}$ do not quite form a simplicial map as $\phi$ respects the algebra structure only up to homotopy, as explained above. 
However, as in the case of the mapping class groups, one can use the explicit homotopies sketched above to rectify $f^i_*$ to a simplicial map, 
while respecting the spectrum structure. 
The rectified map gives the equivalence of spectra.
\end{proof}

Proposition~\ref{infinit1}, Theorem~\ref{LBT} and Theorem~\ref{right} combine to prove Theorem~\ref{main}.

\section{Punctured surfaces}

We have constructed a 2-category $\TT$ which is adequate for proving Theorem~\ref{main}. There are of course many possible versions of $\TT$. 
M.~Weiss suggested the following modification of $\TT$, which has the advantage of being a little more natural to construct, 
but the disadvantage of loosing basepoints, thus making the 
comparison to the original infinite loop space structure on $\Z\x\B\Aut_\infty^+$ more difficult: 

Let $\fD$ be the 2-category obtained from $\TT$ by replacing the 1-morphism $\LL$ by a punctured cylinder $\CC=(S^1\x I)\backslash \{*\}$, with 
one incoming and one outgoing boundary circle. The 1-morphisms of $\fD$ are ---possibly punctured--- cobordisms build out of $\mathbb{P,T,D}$ and $\CC$, 
and the 2-morphisms are homotopy 
equivalences which fix the boundary circles, but not the punctures.

As soon as a surface is punctured, it is homotopy equivalent to a basepoint-free graph $G_{n,k}$, relative to the boundary circles. 
Let $A_{n,k}^0=\pi_0\Htpy(G_{n,k};\tilde\del)$ denote the group of components of the space of homotopy equivalences of $G_{n,k}$ which fix its $k$ boundary 
circles. 
Let $\fD_p(k,1)$ denote the subcategory of $\fD(k,1)$ generated by the punctured surfaces. We have 
$$\D_p(k,1):=\B\fD_p(k,1)\simeq \coprod_{n\ge 0}\B A^0_{n,k+1}\ .$$
As the groups $A_{n,k}^0$ have the same stable homology as $A_{n,k}$ \cite{HatWah04}, one can 
run through Sections 2-3-4 above replacing $\TT$ by $\fD$ and $\LL$ by $\CC$, and show that 
$\Om\B\D\simeq \Z\x\B\Aut_\infty^+$. Only Section 5 does not have a straightforward extension.
 
Let $\M$ denote the mapping class group operad of \cite{Til00}. 
One can define $\M$ in terms of the cobordism category by taking $\M(k)=\ES(k,1)$ and the operad composition induced by composition in $\ES$.
It is shown in \cite{Til00} that $\M$-algebras are infinite loop spaces after group completion. 
The pair of pants multiplication 
does not define a symmetric monoidal structure on $\D_p(0,1)$, but it extends to an action of $\M$. 
One can adapt the proof of Theorem~\ref{right} ---or rather the proof of the main theorem of \cite{Wah04}---
to show the equivalence between the infinite loop space structure of $\Om\B\D$, induced by disjoint union on $\D$,  and that of $\Om\B\D_p(0,1)$, 
induced by the $\M$-algebra structure of $\D_p(0,1)$. 
However, the only way I can see for comparing the infinite loop space structure of $\Om\B\D_p(0,1)$ and the `usual' structure of $\Z\x\B\Aut_\infty^+$
is going through $\Om\B\D$ and $\Om\B\T$ via a middle category with both $\LL$ and $\CC$ as 1-morphisms:\\
$\xymatrix{\Om\B\D_p(0,1) \ar[r] & \Om\B\D \ar[d] &  &  \Om\B\coprod_{n\in \N}\B\Aut(F_n)\\
& \Om\B(\D\cup\T) & \Om\B\T \ar[l] & \ar[l] \Om\B\T_D(0,1) \ar[u] }$.

\vs

Note that $\D_p(0,1)$ is a sub-$\M$-algebra of $\D(0,1)$, and one can show that they have the same group completion, namely $\Z\x\B\Aut_\infty^+$. 
We thus have a map of $\M$-algebras
$$\ES(0,1)\rar \D(0,1)$$
which after  group completion gives another model for the infinite loop map 
$$\Z\x\B\Ga_\infty^+\rar \Z\x\B\Aut_\infty^+.$$

\bibliography{biblio}

\end{document}